\documentclass[10pt,twoside]{article}
\usepackage{graphicx}
\usepackage{amsmath}
\usepackage{Latex-document}

\title{\bf  Non Uniformly Hyperbolic Dynamics: \vskip -2mm H\'enon Maps and Related \vskip -2mm Dynamical Systems\vskip 6mm}

\author{Michael Benedicks\vspace*{-0.5cm}\thanks{Department of Mathematics,
Royal Institute of Technology, S-100\,44 Stockholm, Sweden. E-mail: michaelb@math.kth.se}}

\date{\vspace{-8mm}}

\begin{document}

\maketitle

\thispagestyle{first} \setcounter{page}{255}

\begin{abstract}\vskip 3mm
In the 1960s and 1970s a large part of the theory of
dynamical  systems concerned the case of uniformly hyperbolic or Axiom A
dynamical system and abstract ergodic theory of smooth dynamical
systems. However since around 1980 an emphasize
has been on concrete examples of one-dimensional dynamical systems with abundance
of  chaotic behavior (Collet\,\&Eckmann and Jakobson). New proofs of Jakobson's
one-dimensional results were given by Benedicks and Carleson
\cite{BC85} and were considerably extended to apply to the case of H\'enon maps
by the same authors \cite{BC91}. Since then there has been a
considerable development of these techniques and the methods have been
extended to  the ergodic theory and also to other dynamical systems
(work by Viana, Young, Benedicks and many others). In the cases when
it  applies one can now say that this theory  is now almost as complete as the
Axiom A theory.
\vskip 4.5mm


\end{abstract}

\vskip 12mm

\section{Introduction} \setzero
\vskip-5mm \hspace{5mm}

Dynamical systems as a discipline was born in Henri Poincar\'e's famous
treatise of the three body problem. In retrospect arguably one can
view as his most remarkable discovery of the {\it homoclinic
  phenomenon}. Stable and unstable manifolds of a fixed point or
periodic point may intersect at a {\it homoclinic point} thereby
producing a very complicated dynamic behavior---what we now often
call chaotic.

In my opinion in the development of theory of dynamical
systems---like in the development of all good mathematics---one
can clearly see two stages. The first stage is the understanding
of concrete examples and the second stage is generalization. A
large part of the second encompasses the introduction of the right
concepts which makes the arguments in the concrete examples into a
theory.

One such development starts with the two famous papers by M.~L.~Cartwright and
J.~Littlewood ``On non-linear differential equations of the second
order I and II'', which was among the first to treat nonlinear
differential equations in depth. Littlewood was astonished by the
difficulties that arose in studying these model problems and called
the second of these papers ``the monster''.

S.~Smale gives in ``Finding a horseshoe on the beaches of Rio'' an
entertaining account of how he was lead to his now ubiquitous horseshoe
model for chaotic dynamics. In fact in his first paper in dynamical
systems he made the conjecture that ``chaotic dynamics does
not exist'' but received a letter from N.~Levinson with a paper
clarifying the previous work by Cartwright and Littlewood and
which effectively contained a counterexample to Smale's conjecture.
Levinson's paper contained extensive calculations which Smale found
difficult follow and this lead him to construct a model with
minimal complexity but still with the main features of Levinson's ODE model.

Starting from the horseshoe model, Smale and his group at Berkeley started to develop the
theory of uniformly hyperbolic or Axiom A dynamical systems in the seventies.
One central concept in this theory is that of an axiom A attractor for
a diffeomorphism $f$ of a manifold $M$.

Let $\Lambda$ be an invariant set for a diffeomorphism of a manifold
$M$. $\Lambda$ is said to be a {\it hyperbolic set} for $f$ if there is a
continuous splitting of the tangent bundle of $M$ restricted to
$\Lambda$, $TM|_\Lambda$, which is invariant under the derivative map
$Df$: $TM|_\Lambda=E^s\oplus E^u$; $Df(E^s)=E^s$; $Df(E^u)=E^u$; and for
which there are constants $C>0$ and $c>0$, such that $||Df^n|_{E^s}||
\leq Ce^{-cn}$  and  $||Df^n|_{E^u}||\leq Ce^{-cn}$, for all $n\geq0$ and
there is a uniform lower bound for the angle between
stable and unstable manifolds: $\text{angle}(E^u(x),E^s(x))\geq C,\
\forall x\in \Lambda$.

$\Lambda$ is called an {\it attractor} (in the sense of Conley) if there is a
neighborhood $U \supset \Lambda$ such that $\overline{f(U)}\subset U$
and $\Lambda=\bigcap_{n=1}^\infty f^n(U)$. This attractor is {\it
topologically transitive} if there is $x\in\Lambda$ with dense
orbit in $\Lambda$. If moreover $f|_\Lambda$ is uniformly hyperbolic $\Lambda$ is an
{\it Axiom A attractor}.

The ergodic theory of Axiom A attractors was developed by
Sinai, Ruelle and Bowen in the 1970. In particular for an Axiom A
attractor they constructed so called SRB-measures: measures with
absolutely continuous conditional measures on unstable manifolds. This
measures $\mu$ are also {\it physical measures} in the sense of Ruelle since
for $z_0$ in a set of initial points of positive Riemannian measure
the Birkhoff sums $\frac{1}{n}\sum_{j=0}^{n-1}\delta_{f^jz_0}\to \mu\quad\text{as }n\to\infty$.

In fact for a topologically transitive Axiom A attractor $\mu$ is
unique and the Birkhoff sums converges for a.e. $z_0\in B$, where $B$ is the
{\it basin of attraction} $B=\bigcup_{j\geq 0}f^{-j}(U)$.
Although the Axiom A theory is quite satisfactory and complete it is
not applicable to very many concrete dynamical systems. There is a
general theory of ergodic theory of smooth dynamical systems due to
among others Pesin, Katok, Ledrappier and others but concrete examples
were lacking. Around 1980 the theory of chaotic one-dimensional maps
was started. M.~Misiurewicz studied multimodal maps $f$ of the interval,
whose critical set  or set of turning points, ${\mathcal C}$, has the
property that for all $z_0\in{\mathcal C}$ and all $j\geq1$,
$\text{dist}(f^jz_0,{\mathcal C})\geq\delta>0$ and proved existence of absolutely
continuous invariant measures. Then Collet and Eckmann \cite{CE80}
proved abundance, i.e. positive Lebesgue measure, of aperiodic behavior for
a family of unimodal maps of the interval and, Jakobson in \cite{Ja81}
proved abundance of existence of absolutely continuous invariant measures for the
quadratic family. A new proof of Jakobson's theorem was then given by
Lennart Carleson and myself in \cite{BC85} and the methods from this
paper were later used by us in \cite{BC91} to prove aperiodic, chaotic behavior for
a class of H\'enon maps, which are small perturbations of quadratic maps.
The methods of \cite{BC91} have turned out to be useful for several other
dynamical systems and the corresponding ergodic theory has been developed.
There have been other accounts of this development, in particular by my
collaborators Lai-Sang Young (see e.g. \cite{Yo00}) and Marcelo Viana (see e.g.
the proceedings of ICM98 \cite{Vi98}).

\section{H\'enon maps}

\vskip-5mm \hspace{5mm}

In 1978, M.~H\'enon  proposed as a model for non-linear two-dimensional
dynamical systems the map
$$
(x,y)\mapsto (1+y -ax^2,bx)\qquad 0<a<2,\ b>0.
$$
He chose the parameters $a=1.4$ and $b=0.3$ and proved that
$f=f_{a,b}$ has an attractor in the sense of Conley.

He also verified numerically that this H\'enon map has sensitive
dependence on initial conditions and produced his well-known computer
pictures of the attractor. H\'enon proposed
that this dynamical system should have a ``strange attractor'' and
that it should be more eligible to analysis than the ubiquitous Lorenz system.

In principle most initial points could be attracted to a long periodic
cycle. In view of the famous result of S.~Newhouse, \cite{Ne74},
periodic attractors are topologically generic, so it was not at all a
priori clear that the attractor seen by H\'enon was not a long stable
periodic orbit.

However Lennart Carleson and I, \cite{BC91}, managed to prove that
what H\'enon conjectured was true---not for the parameters
$(a,b)=(1.4,0.3)$ that H\'enon studied---but for small $b>0$. In
fact we managed to prove the following result:

{\bf Theorem 1.} {\it  There is a constant $b_0>0$ such that for
all $b$, $0<b <b_0$ there is a set $A_b$ of parameters $a$,  such
that its one-dimensional Lebesgue measure $|A_b|>0$ and such that
for all $a\in A_b$, $f=f_{a,b}$, has the following properties

\begin{enumerate}
\item  There is an open set $U=U_{a,b}$ such that
$\overline{f(U)} \subset U$ and $\Lambda=\bigcap_{n=0}^\infty
f^n(U)=\overline{W^u(\hat{z})}$, where $W^u(\hat{z})$ is the unstable manifold
of the fixed point $\hat{z}$ of $f$ in the first quadrant.

\item There is a point $z_0=z_0(a,b)$ such that $\{f^j(z_0)\}_{j=0}^\infty$
is dense on $\Lambda$, and there is $c>0$ such that
$|Df^j(z_0)(0,1)|\geq e^{cj},\ j=1,2,\dots$\,.

\end{enumerate}}

Hence $\Lambda$ is a topological transitive attractor with sensitive
dependence on initial conditions. An immediate consequence of Fubinis
theorem is that the ``good parameter set'' $A=\bigcup_{b>0}A_b\times\{b\}$
is (a Cantor set) of positive two-dimensional Lebesgue measure.

The first part of the theorem is  easy to prove and the result is true
for an open set of parameters, i.e. for a small rectangle close to $a=2$ and $b=0$
contained in $\{(a,b):0<a<2,\,b>0\}$. The system is dissipative since
$|\text{det}Df_{a,b}|= |b|<1$. In this case applying an argument of Palis and Takens,
\cite{PT93}, it follows that a region that is enclosed
by pieces of stable and unstable manifolds of the same fixed point $\hat{z}$
is attracted to the unstable manifold $W^u(\hat z)$. With
some additional arguments one can see that also a neighborhood of the
closure of the unstable manifold is attracted. However the attractor
could a priori be a proper subset of $\overline{W^u(\hat{z})}$.

The second part of the theorem is only true for parameters $(a,b)\in A$.
The main ingredient in the proof of the second part of the theorem
is the identification of a {\it critical set} ${\mathcal C}$ for these H\'enon
attractors. The set ${\mathcal C}$ is countable and located on
$W^u(\hat{z})$ but it is natural to expect that the Hausdorff
dimension of $\overline{\mathcal C}$ is positive.

For each $z_0\in {\mathcal C}$ the following holds:
(i) $|Df^j(z_0)\tau(z_0)| \leq e^{-cj}\ \forall j \geq 1$,
where $\tau(z_0)$ is the tangent vector of the unstable manifold at
$z_0$; (ii) trough each $z_0\in {\mathcal C}$ there is a local unstable
  manifold $W^s(z_0)$, tangent to $W^u(\hat{z})$.

The proof of the theorem is a huge induction in time $n$. Successively
preliminary  critical points or sets, {\it precritical points},   are
defined on higher and   higher {\it generations} of the unstable
manifold. We roughly say that  $z\in W^u(\hat{z})$ is of generation
$G$ if $z\in f^G(\gamma)\setminus f^{G-1}(\gamma)$, where $\gamma$ is
the horizontal segment of the   local unstable manifold $W^u(\hat{z})$
through the fixed point.

In analogy with the methods from the one-dimensional case of
\cite{BC85} parameters are chosen so that inductively
$d(f^jz_0,{\mathcal C}_G)\geq e^{-\alpha j}\ \forall j \leq n,\
\forall z_0\in{\mathcal C}_G$, where $C_G$ is the set of precritical points of generation $\leq
G=\theta(b)\cdot n$, where $\theta(b)= C/\log(1/b)$ and
$\alpha>0$ is a suitably chosen numerical constant. Moreover a
further parameter selection is made so that, informally speaking, too
deep returns close to the critical set do not occur too often. The
estimate of the measure of the set, deleted because of this condition, is
made by a large deviation argument.

\section{The ergodic theory} \setzero

\vskip-5mm \hspace{5mm}

{\it Existence of SRB-measures.} For the set of ``good parameters''
$A$ of Theorem 1, Lai-Sang Young and I proved in \cite{BeY92}, the following result

{\bf Theorem 2.} {\it   For all $(a,b)\in A$, $f_{a,b}$ has  a
unique SRB measure supported on the attractor.}

As a consequence it follows by general smooth ergodic theory that
there is a set of initial points $E$ of positive two-dimensional
Lebesgue measure, a subset of the topological basin $B$, such that for all
$z_0\in E$ the Birkhoff sums
$\frac{1}{n}\sum_{j=0}^{n-1}\delta_{f^jz_0} \to \mu$, weak-$*$ as $n\to\infty$.

{\it Decay of correlation and the central limit theorem.} For the same
class of H\'enon maps as in Theorem 1, L.S.~Young and I, \cite{BeY00},
managed to prove decay of correlation and a version of the
central limit theorem. Our main results may be summarized in the
following theorem:

{\bf Theorem 3.} {\it   Suppose $\varphi$ and $\psi$ are H\"older
observables, i.e. they are functions on the plane that belong to
some H\"older class $\alpha,\ 0<\alpha\leq 1$. Then for some $C>0$
and $c>0$
\begin{enumerate}
\item $\left|\int\varphi(f^j(x))\,\psi(x)\,d\mu -\left(\int\varphi\,d\mu\right)
(\int\psi\,d\mu)\right| \leq Ce^{-cj}\qquad \forall j \geq 0$;

\item
  $\mu\left(\{x:\frac{1}{\sqrt{n}}\left(\sum_{j=0}^{n-1}\varphi(f^jx)-n\int\varphi\,
d\mu\right)\leq
    t\}\right)\to \Phi_\sigma(t)$ as $n\to\infty$, where $\Phi_\sigma(t)$ is the
  normal distribution function ${\mathcal N}(0,\sigma)$.
\end{enumerate}}

The methods used to prove this theorem involved the definition of a
return set $X=X_u\cap X_s$, where $X_u$ is a set of approximately
horizontal long unstable manifolds $\gamma_u$ and $X_s$ is a set of approximately
vertical stable manifolds $\gamma_s$, indexed, say, by an arclength coordinate
of its intersection with one of the unstable manifolds of $X_u$. A
dynamical tower construction was made and a Markov extension (Markov
partition on the tower)
was constructed. One of the key estimates concerns the distribution of
the return time
$R(x)=R_i$ for $x\in X^{i}$, which is defined as the time the image a partition element
$X^i$ returns to the base of the tower. $R_i$ may be defined as the
first time such that $f^{R_i}(X^i)\supset (\gamma_u\cap X)$ for all
$\gamma_u$ such that $f^{R_i}(X^i)\cap\gamma_u \neq \emptyset$ and
$R(x)\in X$. Our estimate is that there are constants $C$ and $c>0$
such that for each $\gamma_u$, $|\{x\in X\cap\gamma_u: R(x)>t\}|\leq
C^{-ct}$. L.S.~Young was then able to give a generalization of this setting
of dynamical towers, which applies to other dynamical systems. In
particular she managed to prove exponential decay of correlation for dissipative billiards (\cite{Yo98}).

{\it The metric basin problem.} A natural question that arises in
connection with Theorem 2 is for which set of initial points $z_0$ the
Birkhoff sums $n^{-1}\sum_{j=0}^{n-1} \delta_{f^jz_0}$ converge
weak-$*$ to the SRB-measure $\mu$. As previously mentioned from the
smooth ergodic theory it only follows that this is true for a set of initial
points of positive two-dimensional Lebesgue measure.

In \cite{BeV01}, Marcelo Viana and I were able to complete the picture: in
fact almost all points of the topological basin are generic for
the SRB-measure and the basin is foliated by stable manifolds.

{\bf Theorem 4.} {\it   Let us consider the set of H\'enon maps
$f_{a,b}$, where  $(a,b)\in A$ (the set of good parameters of
Theorem 1). Then the following holds

\begin{enumerate}

\item Through Lebesgue a.e. $z_0\in B$ there is an infinitely long stable manifold
$W^s(z_0)$ that hits the attractor.

\item For a.e. initial point $z_0\in B$,
  $\frac{1}{n}\sum_{j=0}^{n-1}\delta_{f^jz_0}\to \mu$ weak-$*$ as $n\to\infty$,
where $\mu$ is the SRB-measure of Theorem 2.

\end{enumerate}}

This work was in fact carried out in the more general setting of the
H\'enon-like maps of Mora\,\&Viana (see  Section \ref{sec:other} below).
We were also able to characterize the topological basin in this
setting: its boundary is the stable manifold of the fixed point in the
third quadrant (this was proved independently by Y.~Cao \cite{Cao99}).

\section{Other dynamical systems}\label{sec:other}

\vskip-5mm \hspace{5mm}

{\it H\'enon-like maps.} After a rescaling of the second coordinate, the
H\'enon maps may be written as $(x,y)\mapsto (1-ax^2, 0)+\sqrt{b}(y,x)$.
More generally  Mora\,\&Viana, \cite{MV93} con-sided maps of the form
$f_a(x,y)=(1-ax^2, 0)+\psi(x,y)$, where $c_1 b\leq |\text{det}(Df(x,y)|
\leq c_2 b$, $||D(\log|\det Df|)||_{\infty}\leq C$ for some $c_1$,
$c_2$, $C$ and $||\psi||_{C^3}={\mathcal O}(b^{\frac12})$. They managed
to carry out the same program as in Theorem 1 for these class of maps,
which they called {\it H\'enon-like maps}, and managed to prove
prevalence of strange attractors, i.e. that there is a set of positive
measure of parameters $a$ so that $f_a$ exhibits a strange attractor.

Mora\,\&Viana used a one-parameter family of maps $g_\mu$, $-\varepsilon<\mu<\varepsilon$,
such that $g_0$ has a homoclinic tangency and proved that there is a positive
measure set of parameters $\mu$ and a neighborhood $U_\mu$ such that
$g^N_\mu|_{U_\mu}$ has a H\'enon-like  strange attractor. This is done
following Palis\,\&Takens \cite{PT87}, by proving that
there is a linear change of variables and in the  parameters $\Phi_N$ so that
$\Phi_N^{-1}\circ g_\mu^N \circ \Phi_N(\xi,\eta)=(1-a\xi^2, 0)+\psi(\xi,\eta)$,
where $\psi$ satisfies the appropriate estimates of the H\'enon-like
maps. Note that the consequence of this is not the existence of a
global attractor as in Theorem 1 but the existence of a local
attractor close to the homoclinic tangency. (A homoclinic tangency does
really occur close to H\'enon's classical parameters $a=1.4$ and
$b=0.3$: this was proved by Fornaess and Gavasto, \cite{FoG99}.)

{\it Saddle node bifurcations.} Another case where methods based on
those in \cite{BC91} turned out to be useful is the case of saddle
node bifurcations treated by Diaz, Rocha and Viana in
\cite{DRV96}. In that paper they show that when unfolding a one-parameter family with a critical
saddle node cycle H\'enon like strange attractors appear with positive density
at the bifurcation value. Moreover they prove that in an open class of
such families the strange attractors are of global type.

This work was continued by M.J.~Costa \cite{Costa01},
\cite{Costa98b}, who proved that global strange attractors also appear
when destroying a hyperbolic set (horseshoe) by collapsing it with a
periodic attractor.

{\it Viana's dynamical systems with multiple expanding directions}. In
\cite{Vi97}, M.~Viana studied the  dynamical system
$f:{\bf T}\times I\mapsto {\bf T}\times I$ given by the the following
skew product $(\theta, x)\mapsto(m\theta \ (\text{mod}\ 1),a_0(\theta)-x^2)$,
where $a(\theta)=a_0+\alpha\,\sin\pi(\theta-\frac{1}{2})$.

If  $m$ is a sufficiently large integer ($\geq 16$ is enough), $a_0$ is such
that $f_{a_0}(x)=a_0-x^2$ is a Misiurewicz map,
i.e. $|f^j(0)|\geq\delta,\ \forall j\geq 1$, and $\alpha>0$ is
sufficiently small, he managed to prove that for a.e. $(\theta,x)$
and some constants $C$ and $c>0$, $|\partial_x f^j(\theta,x)|\geq Ce^{cj},\ \forall j\geq1$.
In the second part of this paper Viana also considers skew products of
H\'enon maps driven by circle endomorphisms.

An important difference with the situation in this paper compared to
earlier work in the area is that the exponential approach rate
condition of an orbit relative to the critical set is no longer
satisfied. Instead this is replaced by a statistical property: very
deep and very frequent returns to the critical region is very
unlikely. The argument is based on an extension of the large deviation
argument from \cite{BC91}.

SRB measures for the Viana maps were constructed by J.F.~Alves in \cite{Al00}.
An important concept introduced by Alves was the
notion of {\it hyperbolic times}, which are a generalization of the
{\it escape situations} that were considered in \cite{BC91} and they
are also similar to the base in the tower construction in \cite{Yo98}, \cite{BeY00}.

{\it Infinite-modal maps and flows.} Motivated by the study of
unfolding of saddle-focus connections for flows in three dimensions
Pacifico, Rovella and Viana, \cite{PRV98}, studied parameterized families of one-dimensional maps
with infinitely many critical points. They prove that for a positive
Lebesgue measure set of parameter values the map is transitive and
almost all orbits have positive Lyapunov exponent. There has been
a considerable amount of work on flows, by Viana, Luzatto,
Pumari\~no, Rodr\'{\i}guez  and others,
 where techniques from \cite{BC91} have played an important role.
For a survey of these and related results I refer to \cite{Vi94a}.

In a different direction is the recent proof by W.~Tucker, \cite{Tuc2002}, of the
existence of chaotic behaviour for the Lorenz attractor.

{\it  The attractors of Wang and Young.} In a recent paper, D.~Wang and
L.S.~Young, \cite{WY01}, carry out a theory of attractors which generalize
the H\'enon maps in a different direction. They consider maps of a
two-dimensional manifold of the form
$f(x,y,a,b)=(F(x,y,a,b),0)+b\left(u(x,y,a,b),v(x,y,a,b)\right)$.
This class clearly differs from the H\'enon-like maps of
Mora\,\&Viana. In particular the theory applies
to perturbations of certain one-dimensional multimodal maps, with a
transversality condition in the parameter dependence. The techniques
are analogous to these in \cite{BC91} but
more information on the geometric structure, in particular of the
critical set, is achieved. Most
previous results are obtained in this setting but also  new results,
e.g. on some similar dynamical systems and on topological entropy.

\section{Random perturbations and stochastic stability}

\vskip-5mm \hspace{5mm}

A natural
question is how the
statistical properties of a dynamical system with chaotic behavior
behaves when it is perturbed randomly by some small noise at each
iterate. Here we will mainly consider independent additive noise
and assume that the underlying ambient space $M$ is either a subset of
Euclidean space or a torus but cases of more general manifolds and
more general perturbations can also be considered  (for this see
several papers and books by Y.~Kiefer).

Let $f:M\to M$ and  suppose that $\xi_n$, $n\geq 0$, are independent
identically distributed random
variables with an absolutely continuous probability density supported
in a small ball $B(0,\varepsilon)$ around the origin and consider the
Markov chain $\{X_n\}_{n=0}^\infty$ defined by $X_{n+1}=f(X_n)+\xi_n$.
Then there is a stationary transition probability
$p_\varepsilon(E|x)=p(X_{n+1}\in E | X_n=x)$ and also at
least some stationary measure $\nu_{\varepsilon}$ satisfying the
fixed point equation $\nu_\varepsilon(E)=\int p_\varepsilon(E|x)\,d\nu_\varepsilon(x)$.

The obvious questions are here whether $\nu_\varepsilon$ is unique and
in that case if  $\nu_\varepsilon$ tends  to an invariant measure
of the unperturbed system when $\varepsilon\to 0$. This is the problem
of {\it Stochastic Stability}. For the case of Axiom A attractors such
results were proved by Y.~Kiefer and L.S.~Young.

Now  such results have also been obtained for the non-uniformly hyperbolic
dynamical systems described above. In the case of the quadratic interval maps
of \cite{BC85}, L.S.~Young and I proved in \cite{BeY92}, under suitable
assumptions on
the density of the perturbations, that $\nu_\varepsilon$ is unique and
$\nu_\varepsilon\to\mu$
weak-$*$, as $\varepsilon\to0$, where $d\mu=\varphi\,dx$ is the absolutely continuous
invariant measure. V.~Baladi and M.~Viana, \cite{BaV96}, managed to improve this to
prove that the density of $\nu_\varepsilon$, $\psi_\varepsilon$,
converges to $\varphi$ in $L^1$-norm.

M.~Viana and I have recently proved results on weak-$*$ stochastic
stability for the H\'enon maps of \cite{BC91} and the H\'enon-like
maps of \cite{MV93} (see \cite{BeV2}).
For a recent paper  on decay   of correlation for {\it random skew
  products} of quadratic maps see Baladi\,\&Benedicks\,\&Maume-Deschampes \cite{BBM}.

\section{Open problems and concluding remarks}

\vskip-5mm \hspace{5mm}

The most important problem in this general
area is the problem of
positive Lyapunov exponent for the Standard Map, i.e. the map of the
two-dimensional torus defined by
$(x,y)\mapsto(2x-y+K\sin2\pi  x,x) \quad (\text{mod } {\bf Z}^2)$.

The general belief is that at least for some parameters $K$ there is
at least one ergodic component of positive Lebesgue measure with
positive Lyapunov exponent. Nothing is however rigorously known in
spite of intensive work by many people. One of the most interesting
results by Duarte \cite{Du94} goes in the opposite direction: for a
residual set of parameters $K$ the closure of the elliptic points can
have Hausdorff dimension arbitrarily close to 2.

One important difference between the Standard Map and the H\'enon maps for the
good parameters is that ``the critical set'' in the H\'enon case is
rather small. It has a hierarchical structure and conjecturally the
Hausdorff dimension of its closure should be ${\mathcal O}(1/(\log(1/b))$.
The critical set for the standard map (if possible to define)
should have Hausdorff-dimension $\geq 1$.

A.~Baraviera proved in his recent thesis \cite{Bar00}, positive
Lyapunov exponent for the Standard Map with parameters driven by an
expanding circle endomorphism.

On natural class of problem is to consider are {\it skew products} of Viana's
type, where the parameters are driven by more general maps. One
possible choice is to let the driving map be a non-uniformly
hyperbolic quadratic map, either a Misiurewicz map or more generally
the class of quadratic maps of \cite{BC85}, \cite{BC91}. A more
difficult project would be to study the case when driving map is a circle rotation.

{\it The general picture for dissipative H\'enon maps}.
It is a natural question to consider what happens for other parameters
than the ones considered in \cite{BC91}. One possible scenario is
outlined in the following questions, which are much related to
J.~Palis conjectures, \cite{Pa00}, in the $C^r$-generic setting.

{\bf Question 1.} Are there  for Lebesgue almost every parameter
$(a,b)$ in
  $\{(a,b)\in{\mathbf R}^2: 0<a<2,\, b>0\}$  at most finitely many
coexisting strange attractors and stable periodic orbits?

If this is true the Newhouse phenomenon of infinitely coexisting
sinks and Colli's situation \cite{Colli98} of infinitely many coexisting
H\'enon-like strange attractors would only appear for a Lebesgue zero
set of parameters $(a,b)$?

{\bf Question 2.} For the parameter values for which only finitely
many H\'enon-like attractors or sinks coexists: do the respective
basins cover Lebesgue almost all points of the phase space?

{\bf Question 3.} Is the set of parameters for which the H\'enon
map $f_{a,b}$ is hyperbolic dense in the parameter space?

This result if true would correspond to the real Fatou conjecture
proved fairly recently by Graczyk\,\&Swiatek, \cite{GS97}, and Lyubich
\cite{Ly97}. A positive answer to Question 1 would correspond to
Lyubich's result that almost all points in the quadratic family is
either regular or stochastic, \cite{Ly00}.

As can be seen from the above Dynamical Systems is a beautiful mixture
of Topology and Analysis. To the hard analysis of Cartwright, Littlewood
and Levinson, Smale was able to provide a topological
counterpart. Recently again more analytical methods have played an
important role. What will be next?

\label{lastpage}

\end{document}